\documentclass{amsproc}
\usepackage{amsfonts}

\newtheorem{theorem}{Theorem}

\newtheorem{definition}[theorem]{Definition}

\input{tcilatex}

\begin{document}
\title[On the Solutions of Systems of Difference Equations via Tribonacci
Numbers]{On the Solutions of Systems of Difference Equations via Tribonacci
Numbers}
\thanks{}
\author[\.{I}nci Okumu\c{s}, Y\"{u}ksel Soykan]{\.{I}nci Okumu\c{s}, Y\"{u}%
ksel Soykan}
\maketitle

\begin{center}
\textsl{Zonguldak B\"{u}lent Ecevit University, Department of Mathematics, }

\textsl{Art and Science Faculty, 67100, Zonguldak, Turkey }

\textsl{e-mail: inci\_okumus\_90@hotmail.com \ (corresponding author)}

\textsl{\ yuksel\_soykan@hotmail.com}

\textbf{Abstract}
\end{center}

The main objective of this paper is to investigate the explicit form,
stability character and global behavior of solutions of the following two
systems of rational difference equations%
\begin{equation*}
x_{n+1}=\frac{\pm 1}{y_{n}\left( x_{n-1}\pm 1\right) +1},\ \ y_{n+1}=\frac{%
\pm 1}{x_{n}\left( y_{n-1}\pm 1\right) +1},\text{ \ }n=0,1,...
\end{equation*}%
such that their solutions are associated with Tribonacci numbers.

\textbf{2010 Mathematics Subject Classification. }39A10, 39A30

\textbf{Keywords:} \textbf{difference equations, solution, equilibrium
point, tribonacci number, global asymptotic stability.}

\section{Introduction}

Difference equation or discrete dynamical system is a diverse field which
impact almost every branch of pure and applied mathematics. Lately, there
has been great interest in the study of solving difference equations and
systems of difference equations, see [\ref{Yazlik2013}-\ref{Akrour2019}]. In
these studies, the authors deal with the closed-form, stability,periodicity,
boundeness and asymptotically behavior of solutions of nonlinear difference
equations and systems of difference equations. There are many recent
investigations and interest in the field which difference equations have
been studied by several authors, as in the examples given below:

Tollu et al. [\ref{Tollu2014}] considered the following four Riccati
difference equations%
\begin{equation}
x_{n+1}=\frac{1+x_{n}}{x_{n}},\ \ y_{n+1}=\frac{1-y_{n}}{y_{n}}\text{, \ }%
u_{n+1}=\frac{1}{u_{n}+1},\ \ v_{n+1}=\frac{1}{v_{n}-1}\text{,}
\label{eu:fourequ}
\end{equation}%
in which the initial conditions are real numbers. They derived the formulae
for the solutions of equations (\ref{eu:fourequ}) with their solutions are
associated to Fibonacci numbers.

Also, they in [\ref{Tollu2014a}] studied the systems of difference equations%
\begin{equation*}
x_{n+1}=\frac{1+p_{n}}{q_{n}},\ \ y_{n+1}=\frac{1+r_{n}}{s_{n}}\text{, \ }%
n\in 
\mathbb{N}
_{0}\text{,}
\end{equation*}%
where each of the sequences $p_{n}$, $q_{n}$, $r_{n}$ and $s_{n}$ is some of
the sequences $x_{n}$ or $y_{n}$ by their own. They solved fourteen systems
out of sixteen possible systems. In particularly, the representation
formulae of solutions of twelve systems were stated via Fibonacci numbers.

In [\ref{Matsunaga2018}], Matsunaga and Suzuki studied the following system
of rational difference equations%
\begin{equation}
x_{n+1}=\frac{ay_{n}+b}{cy_{n}+d}\text{, \ }y_{n+1}=\frac{ax_{n-1}+b}{%
cx_{n}+d}\text{, \ }n=0,1,...\text{,}  \label{equ:sdllh}
\end{equation}%
where the parameters $a$, $b$, $c$, $d$ and the initial values $x_{0}$, $%
y_{0}$ are real numbers. They obtained the explicit solutions of system (\ref%
{equ:sdllh}) which are generalized Fibonacci sequence.

In [\ref{Ocalan2019}], \"{O}calan and Duman considered the following
nonlinear recursive difference equation%
\begin{equation}
x_{n+1}=\frac{x_{n-1}}{x_{n}},\ n=0,1,...\text{,}  \label{equ:okuhj}
\end{equation}%
with any nonzero initial values $x_{-1}$ and $x_{0}$. Then, they extended
their all results to solutions of the following nonlinear recursive equations%
\begin{equation}
x_{n+1}=\left( \frac{x_{n-1}}{x_{n}}\right) ^{p},\ p>0\text{ and }n=0,1,...%
\text{,}  \label{equ:uhygtf}
\end{equation}%
with any nonzero initial values $x_{-1}$ and $x_{0}$.

Hence, in this study, we consider the following systems of difference
equations%
\begin{equation}
x_{n+1}=\frac{1}{y_{n}\left( x_{n-1}+1\right) +1},\ \ y_{n+1}=\frac{1}{%
x_{n}\left( y_{n-1}+1\right) +1},\text{ \ }n=0,1,...\text{,}
\label{equ:oooooooooooooojk}
\end{equation}%
\begin{equation}
x_{n+1}=\frac{-1}{y_{n}\left( x_{n-1}-1\right) +1},\ \ y_{n+1}=\frac{-1}{%
x_{n}\left( y_{n-1}-1\right) +1},\text{ \ }n=0,1,...\text{,}
\label{equ:waszxb}
\end{equation}%
such that their solutions are associated with Tribonacci numbers.

Our aim in this study is to determine some relationships both between
Tribonacci numbers and and solutions of the aforementioned systems of
difference equations and between the tribonacci constant and the equilibrium
points of these systems of difference equations.

\section{Preliminaries}

\subsection{Linearized stability}

Let us introduce the discrete dynamical system:%
\begin{eqnarray}
x_{n+1} &=&f_{1}\left( x_{n},y_{n}\right) \text{,}  \label{equ:tttttt} \\
y_{n+1} &=&f_{2}\left( x_{n},y_{n}\right) \text{,}  \notag
\end{eqnarray}%
$n\in 
\mathbb{N}
$, where $f_{1}:I_{1}\times I_{2}\rightarrow I_{1}$and $f_{2}:I_{1}\times
I_{2}\rightarrow I_{2}$\ are continuously differentiable functions and $%
I_{1} $, $I_{2}$ are some intervals of real numbers. Also, a solution $%
\{x_{n},y_{n}\}_{n=-k}^{\infty }$ of system (\ref{equ:tttttt}) is uniquely
determined by initial values $\left( x_{0},y_{0}\right) \in I_{1}\times
I_{2} $.

\begin{definition}
An equilibrium point of system (\ref{equ:tttttt}) is a point $\left( 
\overline{x},\overline{y}\right) $ that satisfies%
\begin{eqnarray*}
\overline{x} &=&f_{1}\left( \overline{x},\overline{y}\right) \text{,} \\
\overline{y} &=&f_{2}\left( \overline{x},\overline{y}\right) \text{.}
\end{eqnarray*}
\end{definition}

Together with system (\ref{equ:tttttt}), if we consider the associated
vector map%
\begin{equation*}
K=\left( f_{1},x_{n},f_{2},y_{n}\right) \text{,}
\end{equation*}%
then the point $\left( \overline{x},\overline{y}\right) $ is also called a
fixed point of the vector map $K$.

\begin{definition}
Let $\left( \overline{x},\overline{y}\right) $ be an equilibrium point of
the map $K$ where $f_{1}$,and $f_{2}$ are continuously differentiable
functions at $\left( \overline{x},\overline{y}\right) $. The linearized
system of system (\ref{equ:tttttt}) about the equilibrium point $\left( 
\overline{x},\overline{y}\right) $ is%
\begin{equation}
X_{n+1}=K\left( X_{n}\right) =BX_{n}\text{,}  \label{equ:iiiiiiiiii}
\end{equation}%
where%
\begin{equation*}
X_{n}=\left( x_{n},y_{n}\right) ^{T}
\end{equation*}%
and $B$ is a Jacobian matrix of system (\ref{equ:tttttt}) about the
equilibrium point $\left( \overline{x},\overline{y}\right) $.
\end{definition}

The linearized system, associated to system (\ref{equ:tttttt}), about the
equilibrium point $\left( \overline{x},\overline{y}\right) $ is given by%
\begin{equation*}
\left( 
\begin{array}{c}
x_{n+1} \\ 
y_{n+1}%
\end{array}%
\right) =\left( 
\begin{array}{cc}
\frac{\partial f_{1}}{\partial x_{n}}\left( \overline{x},\overline{y}\right)
& \frac{\partial f_{1}}{\partial y_{n}}\left( \overline{x},\overline{y}%
\right) \\ 
\frac{\partial f_{2}}{\partial x_{n}}\left( \overline{x},\overline{y}\right)
& \frac{\partial f_{2}}{\partial y_{n}}\left( \overline{x},\overline{y}%
\right)%
\end{array}%
\right) \left( 
\begin{array}{c}
x_{n} \\ 
y_{n}%
\end{array}%
\right) \text{.}
\end{equation*}

\begin{definition}
Let $\left( \overline{x},\overline{y}\right) $ be an equilibrium point of
system (\ref{equ:tttttt}).

\begin{description}
\item[(a)] An equilibrium point $\left( \overline{x},\overline{y}\right) $
is called stable if, for every $\varepsilon >0$; there exists $\delta >0$
such that for every initial value $\left( x_{0},y_{0}\right) \in I_{1}\times
I_{2}$, with%
\begin{equation*}
\left\vert x_{0}-\overline{x}\right\vert <\delta ,\text{ }\left\vert y_{0}-%
\overline{y}\right\vert <\delta \text{,}
\end{equation*}%
implying $\left\vert x_{n}-\overline{x}\right\vert <\varepsilon $, $%
\left\vert y_{n}-\overline{y}\right\vert <\varepsilon $, for $n\in 
\mathbb{N}
$.

\item[(b)] If an equilibrium point $\left( \overline{x},\overline{y}\right) $
of system (\ref{equ:tttttt}) is called unstable if it is not stable.

\item[(c)] An equilibrium point $\left( \overline{x},\overline{y}\right) $
of system (\ref{equ:tttttt})\ is called locally asymptotically stable if, it
is stable, and if in addition there exists $\gamma >0$ such that%
\begin{equation*}
\left\vert x_{0}-\overline{x}\right\vert <\gamma ,\text{ }\left\vert y_{0}-%
\overline{y}\right\vert <\gamma \text{,}
\end{equation*}%
and $\left( x_{n},y_{n}\right) \rightarrow \left( \overline{x},\overline{y}%
\right) $ as $n\rightarrow \infty $.

\item[(d)] An equilibrium point $\left( \overline{x},\overline{y}\right) $
of system (\ref{equ:tttttt})\ is called a global attractor if $\left(
x_{n},y_{n}\right) \rightarrow \left( \overline{x},\overline{y}\right) $ as $%
n\rightarrow \infty $.

\item[(e)] An equilibrium point $\left( \overline{x},\overline{y}\right) $
of system (\ref{equ:tttttt})\ is called globally asymptotically stable if it
is stable, and a global attractor.
\end{description}
\end{definition}

\begin{theorem}[The Linearized Stability Theorem]
\label{theorem:ertgbns}

Assume that%
\begin{equation*}
X_{n+1}=K\left( X_{n}\right) ,n=0,1,...,
\end{equation*}%
be a systemdifference equations such that $\overline{X}$ is a fixed point of 
$F$.

\begin{description}
\item[(a)] If all eigenvalues of the Jacobian matrix $B$ about $\overline{X}$
lie inside the open unit disk $\left\vert \lambda \right\vert <1$, that is,
if all of them have absolute value less than one, then $\overline{X}$\ is
locally asymptotically stable.

\item[(b)] If at least one of them has a modulus greater than one, then $%
\overline{X}$\ is unstable.
\end{description}
\end{theorem}

\subsection{Tribonacci numbers}

Now, we give information about Tribonacci numbers that we afterwards need in
the paper.

The Tribonacci sequence $\left\{ T_{n}\right\} _{n=0}^{\infty }$ is defined
by the third-order recurrence relations%
\begin{equation}
T_{n+3}=T_{n+2}+T_{n+1}+T_{n}\text{,}  \label{equ:opopopop}
\end{equation}%
with initial conditions $T_{0}=0$, $T_{1}=1$, $T_{2}=1$. Also, it can be
extended the Tribonacci sequence backward (negative subscripts) as%
\begin{equation}
T_{-n}=T_{-n+3}-T_{-n+2}-T_{-n+1}\text{.}  \label{equ:otuotu}
\end{equation}%
It can be clearly obtained that the characteristic equation of (\ref%
{equ:opopopop}) has the form%
\begin{equation}
x^{3}-x^{2}-x-1=0  \label{equ:eeeeeeeeee}
\end{equation}%
such that the roots%
\begin{eqnarray*}
\alpha &=&\frac{1+\sqrt[3]{19+3\sqrt{33}}+\sqrt[3]{19-3\sqrt{33}}}{3} \\
\beta &=&\frac{1+\omega \sqrt[3]{19+3\sqrt{33}}+\omega ^{2}\sqrt[3]{19-3%
\sqrt{33}}}{3} \\
\gamma &=&\frac{1+\omega ^{2}\sqrt[3]{19+3\sqrt{33}}+\omega \sqrt[3]{19-3%
\sqrt{33}}}{3}
\end{eqnarray*}%
where $\alpha $ is called Tribonacci constant and%
\begin{equation*}
\omega =\frac{-1+i\sqrt{3}}{2}=\exp \left( 2\pi i/3\right)
\end{equation*}%
is a primitive cube root of unity. Therefore, Tribonacci sequence can be
expressed using Binet formula%
\begin{equation*}
T_{n}=\frac{\alpha ^{n+1}}{\left( \alpha -\beta \right) \left( \alpha
-\gamma \right) }+\frac{\beta ^{n+1}}{\left( \beta -\alpha \right) \left(
\beta -\gamma \right) }+\frac{\gamma ^{n+1}}{\left( \gamma -\alpha \right)
\left( \gamma -\beta \right) }\text{.}
\end{equation*}

Furthermore, there exist the following limit%
\begin{equation}
\lim_{n\rightarrow \infty }\frac{T_{n+r}}{T_{n}}=\alpha ^{r}\text{,}
\label{lim:wwswws}
\end{equation}%
where $r\in 
\mathbb{Z}
$ and $T_{n}$ is the $n$th Tribonacci number.

\section{Main Results}

In this section, we introduce our results.

\subsection{The System (\protect\ref{equ:oooooooooooooojk})}

In this subsection, we present our main results related to the system (\ref%
{equ:oooooooooooooojk}). Our aim is to investigate the general solution in
exact form of system (\ref{equ:oooooooooooooojk}).and the asymptotic
behavior of solutions of system (\ref{equ:oooooooooooooojk}).

\begin{theorem}
\label{theorem:resxbnr}Let $\left \{ x_{n},y_{n}\right \} _{n=-1}^{\infty }$
be a solution system (\ref{equ:oooooooooooooojk}). Then, for $n=0,1,2,...$,
the form of solutions $\left \{ x_{n},y_{n}\right \} _{n=-1}^{\infty }$ is
given by%
\begin{eqnarray*}
x_{2n-1} &=&\frac{T_{2n-2}x_{-1}y_{0}+\left( T_{2n}-T_{2n-1}\right)
y_{0}+T_{2n-1}}{T_{2n-1}x_{-1}y_{0}+\left( T_{2n-2}+T_{2n-1}\right)
y_{0}+T_{2n}}\text{,} \\
x_{2n} &=&\frac{T_{2n-1}y_{-1}x_{0}+\left( T_{2n+1}-T_{2n}\right)
x_{0}+T_{2n}}{T_{2n}y_{-1}x_{0}+\left( T_{2n-1}+T_{2n}\right) x_{0}+T_{2n+1}}%
\text{,} \\
y_{2n-1} &=&\frac{T_{2n-2}y_{-1}x_{0}+\left( T_{2n}-T_{2n-1}\right)
x_{0}+T_{2n-1}}{T_{2n-1}y_{-1}x_{0}+\left( T_{2n-2}+T_{2n-1}\right)
x_{0}+T_{2n}}\text{,} \\
y_{2n} &=&\frac{T_{2n-1}x_{-1}y_{0}+\left( T_{2n+1}-T_{2n}\right)
y_{0}+T_{2n}}{T_{2n}x_{-1}y_{0}+\left( T_{2n-1}+T_{2n}\right) y_{0}+T_{2n+1}}%
\text{,}
\end{eqnarray*}%
where $T_{n}$ is the $n$th Tribonacci number and the initial conditions $%
x_{-1}$, $y_{-1}$, $x_{0}$, $y_{0}\in 
\mathbb{R}
-F$, with $F$ is the forbidden set of system (\ref{equ:oooooooooooooojk})
given by%
\begin{equation*}
F=\dbigcup \limits_{n=-1}^{\infty }\left \{ \left(
x_{-1},y_{-1},x_{0},y_{0}\right) :A_{n}=0\text{, }B_{n}=0\text{, }C_{n}=0%
\text{, }D_{n}=0\right \}
\end{equation*}%
where%
\begin{eqnarray*}
A_{n} &=&T_{2n-1}x_{-1}y_{0}+\left( T_{2n-2}+T_{2n-1}\right) y_{0}+T_{2n}%
\text{,} \\
B_{n} &=&T_{2n}y_{-1}x_{0}+\left( T_{2n-1}+T_{2n}\right) x_{0}+T_{2n+1}\text{%
,} \\
C_{n} &=&T_{2n-1}y_{-1}x_{0}+\left( T_{2n-2}+T_{2n-1}\right) x_{0}+T_{2n}%
\text{,} \\
D_{n} &=&T_{2n}x_{-1}y_{0}+\left( T_{2n-1}+T_{2n}\right) y_{0}+T_{2n+1}\text{%
.}
\end{eqnarray*}
\end{theorem}

\textbf{Proof.} We use the induction on $k$. For $k=0$, the result holds.
Suppose that $k>0$ and that our assumption holds for $k-1$. That is,%
\begin{eqnarray*}
x_{2k-3} &=&\frac{T_{2k-4}x_{-1}y_{0}+\left( T_{2k-2}-T_{2k-3}\right)
y_{0}+T_{2k-3}}{T_{2k-3}x_{-1}y_{0}+\left( T_{2k-4}+T_{2k-3}\right)
y_{0}+T_{2k-2}}\text{,} \\
x_{k-2} &=&\frac{T_{2k-3}y_{-1}x_{0}+\left( T_{2k-1}-T_{2k-2}\right)
x_{0}+T_{2k-2}}{T_{2k-2}y_{-1}x_{0}+\left( T_{2k-3}+T_{2k-2}\right)
x_{0}+T_{2k-1}}\text{,} \\
y_{2k-3} &=&\frac{T_{2k-4}y_{-1}x_{0}+\left( T_{2k-2}-T_{2k-3}\right)
x_{0}+T_{2k-3}}{T_{2k-3}y_{-1}x_{0}+\left( T_{2k-4}+T_{2k-3}\right)
x_{0}+T_{2k-2}}\text{,} \\
y_{2k-2} &=&\frac{T_{2k-3}x_{-1}y_{0}+\left( T_{2k-1}-T_{2k-2}\right)
y_{0}+T_{2k-2}}{T_{2k-2}x_{-1}y_{0}+\left( T_{2k-3}+T_{2k-2}\right)
y_{0}+T_{2k-1}}\text{.}
\end{eqnarray*}%
From system (\ref{equ:oooooooooooooojk}) and (\ref{equ:opopopop}), it
follows that%
\begin{eqnarray*}
x_{2k-1} &=&\frac{1}{y_{2k-2}\left( x_{2k-3}+1\right) +1} \\
&=&\frac{1}{\frac{T_{2k-3}x_{-1}y_{0}+\left( T_{2k-1}-T_{2k-2}\right)
y_{0}+T_{2k-2}}{T_{2k-2}x_{-1}y_{0}+\left( T_{2k-3}+T_{2k-2}\right)
y_{0}+T_{2k-1}}\left( \frac{T_{2k-4}x_{-1}y_{0}+\left(
T_{2k-2}-T_{2k-3}\right) y_{0}+T_{2k-3}}{T_{2k-3}x_{-1}y_{0}+\left(
T_{2k-4}+T_{2k-3}\right) y_{0}+T_{2k-2}}+1\right) +1} \\
&=&\frac{T_{2k-2}x_{-1}y_{0}+\left( T_{2k-3}+T_{2k-2}\right) y_{0}+T_{2k-1}}{%
\left( T_{2k-4}+T_{2k-3}+T_{2k-2}\right) x_{-1}y_{0}+\left(
T_{2k-2}+T_{2k-1}\right) y_{0}+T_{2k-3}+T_{2k-2}+T_{2k-1}}\text{.}
\end{eqnarray*}%
Therefore, we have%
\begin{equation*}
x_{2k-1}=\frac{T_{2k-2}x_{-1}y_{0}+\left( T_{2k}-T_{2k-1}\right)
y_{0}+T_{2k-1}}{T_{2k-1}x_{-1}y_{0}+\left( T_{2k-2}+T_{2k-1}\right)
y_{0}+T_{2k}}\text{.}
\end{equation*}%
And also, it follows that%
\begin{eqnarray*}
y_{2k-1} &=&\frac{1}{x_{2k-2}\left( y_{2k-3}+1\right) +1} \\
&=&\frac{1}{\frac{T_{2k-3}y_{-1}x_{0}+\left( T_{2k-1}-T_{2k-2}\right)
x_{0}+T_{2k-2}}{T_{2k-2}y_{-1}x_{0}+\left( T_{2k-3}+T_{2k-2}\right)
x_{0}+T_{2k-1}}\left( \frac{T_{2k-4}y_{-1}x_{0}+\left(
T_{2k-2}-T_{2k-3}\right) x_{0}+T_{2k-3}}{T_{2k-3}y_{-1}x_{0}+\left(
T_{2k-4}+T_{2k-3}\right) x_{0}+T_{2k-2}}+1\right) +1} \\
&=&\frac{T_{2k-2}y_{-1}x_{0}+\left( T_{2k-3}+T_{2k-2}\right) x_{0}+T_{2k-1}}{%
\left( T_{2k-4}+T_{2k-3}+T_{2k-2}\right) y_{-1}x_{0}+\left(
T_{2k-2}+T_{2k-1}\right) x_{0}+T_{2k-3}+T_{2k-2}+T_{2k-1}}\text{.}
\end{eqnarray*}%
So, we obtain%
\begin{equation*}
y_{2k-1}=\frac{T_{2k-2}y_{-1}x_{0}+\left( T_{2k}-T_{2k-1}\right)
x_{0}+T_{2k-1}}{T_{2k-1}y_{-1}x_{0}+\left( T_{2k-2}+T_{2k-1}\right)
x_{0}+T_{2k}}\text{.}
\end{equation*}%
Similarly, from system (\ref{equ:oooooooooooooojk}) and (\ref{equ:opopopop}%
), it follows that%
\begin{eqnarray*}
x_{2k} &=&\frac{1}{y_{2k-1}\left( x_{2k-2}+1\right) +1} \\
&=&\frac{1}{\frac{T_{2k-2}y_{-1}x_{0}+\left( T_{2k}-T_{2k-1}\right)
x_{0}+T_{2k-1}}{T_{2k-1}y_{-1}x_{0}+\left( T_{2k-2}+T_{2k-1}\right)
x_{0}+T_{2k}}\left( \frac{T_{2k-3}y_{-1}x_{0}+\left(
T_{2k-1}-T_{2k-2}\right) x_{0}+T_{2k-2}}{T_{2k-2}y_{-1}x_{0}+\left(
T_{2k-3}+T_{2k-2}\right) x_{0}+T_{2k-1}}+1\right) +1} \\
&=&\frac{T_{2k-1}y_{-1}x_{0}+\left( T_{2k-2}+T_{2k-1}\right) x_{0}+T_{2k}}{%
\left( T_{2k-3}+T_{2k-2}+T_{2k-1}\right) y_{-1}x_{0}+\left(
T_{2k-1}+T_{2k}\right) x_{0}+T_{2k-2}+T_{2k-1}+T_{2k}}\text{.}
\end{eqnarray*}%
Thus, we get%
\begin{equation*}
x_{2k}=\frac{T_{2k-1}y_{-1}x_{0}+\left( T_{2k-2}+T_{2k-1}\right) x_{0}+T_{2k}%
}{T_{2k}y_{-1}x_{0}+\left( T_{2k-1}+T_{2k}\right) x_{0}+T_{2k+1}}\text{.}
\end{equation*}%
And also, it follows that%
\begin{eqnarray*}
y_{2k} &=&\frac{1}{x_{2k-1}\left( y_{2k-2}+1\right) +1} \\
&=&\frac{1}{\frac{T_{2k-2}x_{-1}y_{0}+\left( T_{2k}-T_{2k-1}\right)
y_{0}+T_{2k-1}}{T_{2k-1}x_{-1}y_{0}+\left( T_{2k-2}+T_{2k-1}\right)
y_{0}+T_{2k}}\left( \frac{T_{2k-3}x_{-1}y_{0}+\left(
T_{2k-1}-T_{2k-2}\right) y_{0}+T_{2k-2}}{T_{2k-2}x_{-1}y_{0}+\left(
T_{2k-3}+T_{2k-2}\right) y_{0}+T_{2k-1}}+1\right) +1} \\
&=&\frac{T_{2k-1}x_{-1}y_{0}+\left( T_{2k-2}+T_{2k-1}\right) y_{0}+T_{2k}}{%
\left( T_{2k-3}+T_{2k-2}+T_{2k-1}\right) x_{-1}y_{0}+\left(
T_{2k-1}+T_{2k}\right) y_{0}+T_{2k-2}+T_{2k-1}+T_{2k}}\text{.}
\end{eqnarray*}%
Herefrom, we have%
\begin{equation*}
y_{2k}=\frac{T_{2k-1}x_{-1}y_{0}+\left( T_{2k-2}+T_{2k-1}\right) y_{0}+T_{2k}%
}{T_{2k}x_{-1}y_{0}+\left( T_{2k-1}+T_{2k}\right) y_{0}+T_{2k+1}}\text{.}
\end{equation*}

\begin{theorem}
\label{theorem:qqasefrtv}The system (\ref{equ:oooooooooooooojk}) has unique
positive equilibrium point $\left( \overline{x},\overline{y}\right) =\left(
a,a\right) $ and $\left( a,a\right) $ is locally asymptotically stable.
\end{theorem}

\textbf{Proof.} Clearly, equilibrium point of system (\ref%
{equ:oooooooooooooojk}) is the real roots of the equations%
\begin{equation}
\overline{x}=\frac{1}{\overline{x}\left( \overline{y}+1\right) +1}\text{, \ }%
\overline{y}=\frac{1}{\overline{y}\left( \overline{x}+1\right) +1}\text{.}
\label{equ:yykjhgf}
\end{equation}%
In (\ref{equ:yykjhgf}), after some operations, we obtain%
\begin{equation*}
\overline{x}=\overline{y}\text{.}
\end{equation*}%
As a result, we obtain the following equation%
\begin{equation}
\overline{x}^{3}+\overline{x}^{2}+\overline{x}-1=0\text{.}
\label{equ:kkmnbrq}
\end{equation}%
Then, the roots of the cubic equation (\ref{equ:kkmnbrq}) are given by%
\begin{eqnarray*}
a &=&\frac{-1+\sqrt[3]{3\sqrt{33}+17}-\sqrt[3]{3\sqrt{33}-17}}{3}\text{,} \\
b &=&\frac{-1+\omega \sqrt[3]{3\sqrt{33}+17}-\omega ^{2}\sqrt[3]{3\sqrt{33}%
-17}}{3}\text{,} \\
c &=&\frac{-1+\omega ^{2}\sqrt[3]{3\sqrt{33}+17}-\omega \sqrt[3]{3\sqrt{33}%
-17}}{3}\text{,}
\end{eqnarray*}%
where%
\begin{equation*}
\omega =\frac{-1+i\sqrt{3}}{2}=\exp \left( 2\pi i/3\right)
\end{equation*}%
is a primitive cube root of unity.\ So, the root $a$ is only real number.
Therefore, the unique positive equilibrium point of system (\ref%
{equ:oooooooooooooojk}) is $\left( \overline{x},\overline{y}\right) =\left(
a,a\right) $.

Now, we show that the unique positive equilibrium point of system (\ref%
{equ:oooooooooooooojk}) is locally asymptotically stable.

Let $I$ and $J$ are some intervals of real numbers.and consider the functions%
\begin{equation*}
f:I^{2}\times J^{2}\rightarrow I\text{ and }g:I^{2}\times J^{2}\rightarrow J
\end{equation*}%
defined by%
\begin{equation*}
f\left( x_{n},x_{n-1},y_{n},y_{n-1}\right) =\frac{1}{y_{n}\left(
x_{n-1}+1\right) +1},\ g\left( x_{n},x_{n-1},y_{n},y_{n-1}\right) =\frac{1}{%
x_{n}\left( y_{n-1}+1\right) +1}\text{.}
\end{equation*}

We consider the following transformation to build corresponding linearized
form of system (\ref{equ:oooooooooooooojk})%
\begin{equation*}
\left( x_{n},x_{n-1},y_{n},y_{n-1}\right) \rightarrow \left(
f,f_{1},g,g_{1}\right) \text{,}
\end{equation*}%
where%
\begin{eqnarray*}
f\left( x_{n},x_{n-1},y_{n},y_{n-1}\right) &=&\frac{1}{y_{n}\left(
x_{n-1}+1\right) +1}\text{,} \\
f_{1}\left( x_{n},x_{n-1},y_{n},y_{n-1}\right) &=&x_{n}\text{,} \\
g\left( x_{n},x_{n-1},y_{n},y_{n-1}\right) &=&\frac{1}{x_{n}\left(
y_{n-1}+1\right) +1}\text{,} \\
g_{1}\left( x_{n},x_{n-1},y_{n},y_{n-1}\right) &=&y_{n}\text{.}
\end{eqnarray*}

Then, the linearized system of system (\ref{equ:oooooooooooooojk}) about the
equilibrium point $\left( a,a\right) $ under the above transformation is
given as%
\begin{equation*}
X_{n+1}=BX_{n}\text{,}
\end{equation*}%
where $X_{n}=\left( x_{n},x_{n-1},y_{n},y_{n-1}\right) ^{T}$ and $B$ is a
Jacobian matrix of system (\ref{equ:oooooooooooooojk}) about the equilibrium
point $\left( a,a\right) $ and given by%
\begin{eqnarray*}
B &=&\left( 
\begin{array}{cccc}
0 & \frac{-a}{\left( a\left( a+1\right) +1\right) ^{2}} & \frac{-\left(
1+a\right) }{\left( a\left( a+1\right) +1\right) ^{2}} & 0 \\ 
1 & 0 & 0 & 0 \\ 
\frac{-\left( 1+a\right) }{\left( a\left( a+1\right) +1\right) ^{2}} & 0 & 0
& \frac{-a}{\left( a\left( a+1\right) +1\right) ^{2}} \\ 
0 & 0 & 1 & 0%
\end{array}%
\right) \\
&=&\left( 
\begin{array}{cccc}
0 & -a^{3} & a-1 & 0 \\ 
1 & 0 & 0 & 0 \\ 
a-1 & 0 & 0 & -a^{3} \\ 
0 & 0 & 1 & 0%
\end{array}%
\right) \text{.}
\end{eqnarray*}%
Thus, we obtain the characteristic equation of the Jacobian matrix $B$ as%
\begin{equation*}
\left( a^{3}+\lambda ^{2}\right) ^{2}-\left( a-1\right) ^{2}\lambda ^{2}=0%
\text{,}
\end{equation*}%
or%
\begin{equation*}
\left( \lambda ^{2}+\left( a-1\right) \lambda +a^{3}\right) \left( \lambda
^{2}-\left( a-1\right) \lambda +a^{3}\right) =0\text{.}
\end{equation*}%
Hence, it is clearly seen that numerically%
\begin{equation*}
\left \vert \lambda _{1}\right \vert =\left \vert \lambda _{2}\right \vert
=\left \vert \lambda _{3}\right \vert =\left \vert \lambda _{4}\right \vert
=0.40089<1\text{.}
\end{equation*}%
Consequently, the equilibrium point $\left( a,a\right) $ is locally
asymptotically stable. So, this completes the proof.

\begin{theorem}
The equilibrium point of system (\ref{equ:oooooooooooooojk}) is globally
asymptotically stable.
\end{theorem}

\textbf{Proof.} Let $\left \{ x_{n},y_{n}\right \} _{n\geq -1}$ be a
solution system (\ref{equ:oooooooooooooojk}). By Theorem (\ref%
{theorem:qqasefrtv}), we need only to prove that the equilibrium point $%
\left( a,a\right) $ is global attractor, that is%
\begin{equation*}
\lim_{n\rightarrow \infty }\left( x_{n},y_{n}\right) =\left( a,a\right) 
\text{.}
\end{equation*}%
From Theorem (\ref{theorem:resxbnr}), (\ref{equ:eeeeeeeeee}) and (\ref%
{lim:wwswws}), it follows that%
\begin{eqnarray*}
\lim_{n\rightarrow \infty }x_{2n-1} &=&\lim_{n\rightarrow \infty }\frac{%
T_{2n-2}x_{-1}y_{0}+\left( T_{2n}-T_{2n-1}\right) y_{0}+T_{2n-1}}{%
T_{2n-1}x_{-1}y_{0}+\left( T_{2n-2}+T_{2n-1}\right) y_{0}+T_{2n}} \\
&=&\lim_{n\rightarrow \infty }\frac{T_{2n-2}\left( x_{-1}y_{0}+\left( \frac{%
T_{2n}}{T_{2n-2}}-\frac{T_{2n-1}}{T_{2n-2}}\right) y_{0}+\frac{T_{2n-1}}{%
T_{2n-2}}\right) }{T_{2n-1}\left( x_{-1}y_{0}+\left( \frac{T_{2n-2}}{T_{2n-1}%
}+1\right) y_{0}+\frac{T_{2n}}{T_{2n-1}}\right) } \\
&=&\left( \frac{x_{-1}y_{0}+\left( \alpha ^{2}-\alpha \right) y_{0}+\alpha }{%
x_{-1}y_{0}+\left( \frac{1}{\alpha }+1\right) y_{0}+\alpha }\right)
\lim_{n\rightarrow \infty }\frac{T_{2n-2}}{T_{2n-1}} \\
&=&\lim_{n\rightarrow \infty }\frac{T_{2n-2}}{T_{2n-1}} \\
&=&\frac{1}{\alpha } \\
&=&a,
\end{eqnarray*}%
and%
\begin{eqnarray*}
\lim_{n\rightarrow \infty }x_{2n} &=&\lim_{n\rightarrow \infty }\frac{%
T_{2n-1}y_{-1}x_{0}+\left( T_{2n+1}-T_{2n}\right) x_{0}+T_{2n}}{%
T_{2n}y_{-1}x_{0}+\left( T_{2n-1}+T_{2n}\right) x_{0}+T_{2n+1}} \\
&=&\lim_{n\rightarrow \infty }\frac{T_{2n-1}\left( y_{-1}x_{0}+\left( \frac{%
T_{2n+1}}{T_{2n-1}}-\frac{T_{2n}}{T_{2n-1}}\right) x_{0}+\frac{T_{2n}}{%
T_{2n-1}}\right) }{T_{2n}\left( y_{-1}x_{0}+\left( \frac{T_{2n-1}}{T_{2n}}%
+1\right) x_{0}+\frac{T_{2n+1}}{T_{2n}}\right) } \\
&=&\left( \frac{y_{-1}x_{0}+\left( \alpha ^{2}-\alpha \right) x_{0}+\alpha }{%
y_{-1}x_{0}+\left( \frac{1}{\alpha }+1\right) x_{0}+\alpha }\right)
\lim_{n\rightarrow \infty }\frac{T_{2n-1}}{T_{2n}} \\
&=&\lim_{n\rightarrow \infty }\frac{T_{2n-1}}{T_{2n}} \\
&=&\frac{1}{\alpha } \\
&=&a\text{.}
\end{eqnarray*}%
Then, we have%
\begin{equation*}
\lim_{n\rightarrow \infty }x_{n}=a\text{.}
\end{equation*}%
Similarly, we obtain%
\begin{equation*}
\lim_{n\rightarrow \infty }y_{n}=a\text{.}
\end{equation*}%
Therefore, we get%
\begin{equation*}
\lim_{n\rightarrow \infty }\left( x_{n},y_{n}\right) =\left( a,a\right) 
\text{.}
\end{equation*}%
The proof is completed.

\subsection{The System (\protect\ref{equ:waszxb})}

In this subsection, we introduce our main results related to the system (\ref%
{equ:waszxb}). Our aim is to investigate the general solution in explicit
form of system (\ref{equ:waszxb}).and the asymptotic behavior of solutions
of system (\ref{equ:waszxb}).

\begin{theorem}
\label{theorem:kmuyhn}Let $\left \{ x_{n},y_{n}\right \} _{n=-1}^{\infty }$
be a solution system (\ref{equ:waszxb}). Then, for $n=0,1,2,...$, the form
of solutions $\left \{ x_{n},y_{n}\right \} _{n=-1}^{\infty }$ is given by%
\begin{eqnarray*}
x_{2n-1} &=&\frac{-\left( T_{2n-2}x_{-1}y_{0}+\left( T_{2n-1}-T_{2n}\right)
y_{0}+T_{2n-1}\right) }{T_{2n-1}x_{-1}y_{0}-\left( T_{2n-2}+T_{2n-1}\right)
y_{0}+T_{2n}}\text{,} \\
x_{2n} &=&\frac{-\left( T_{2n-1}y_{-1}x_{0}+\left( T_{2n}-T_{2n+1}\right)
x_{0}+T_{2n}\right) }{T_{2n}y_{-1}x_{0}-\left( T_{2n-1}+T_{2n}\right)
x_{0}+T_{2n+1}}\text{,} \\
y_{2n-1} &=&\frac{-\left( T_{2n-2}y_{-1}x_{0}+\left( T_{2n-1}-T_{2n}\right)
x_{0}+T_{2n-1}\right) }{T_{2n-1}y_{-1}x_{0}-\left( T_{2n-2}+T_{2n-1}\right)
x_{0}+T_{2n}}\text{,} \\
y_{2n-1} &=&\frac{-\left( T_{2n-1}x_{-1}y_{0}+\left( T_{2n}-T_{2n+1}\right)
y_{0}+T_{2n}\right) }{T_{2n}x_{-1}y_{0}-\left( T_{2n-1}+T_{2n}\right)
y_{0}+T_{2n+1}}
\end{eqnarray*}%
where initial conditions $x_{-1}$, $y_{-1}$, $x_{0}$, $y_{0}\in 
\mathbb{R}
-F$, with $F$ is the forbidden set of system (\ref{equ:waszxb}) given by%
\begin{equation*}
F=\dbigcup \limits_{n=-1}^{\infty }\left \{ \left(
x_{-1},y_{-1},x_{0},y_{0}\right) :A_{n}=0\text{, }B_{n}=0\text{, }C_{n}=0%
\text{, }D_{n}=0\right \}
\end{equation*}%
where%
\begin{eqnarray*}
A_{n} &=&T_{2n-1}x_{-1}y_{0}-\left( T_{2n-2}+T_{2n-1}\right) y_{0}+T_{2n}%
\text{,} \\
B_{n} &=&T_{2n}y_{-1}x_{0}-\left( T_{2n-1}+T_{2n}\right) x_{0}+T_{2n+1}\text{%
,} \\
C_{n} &=&T_{2n-1}y_{-1}x_{0}-\left( T_{2n-2}+T_{2n-1}\right) x_{0}+T_{2n}%
\text{,} \\
D_{n} &=&T_{2n}x_{-1}y_{0}-\left( T_{2n-1}+T_{2n}\right) y_{0}+T_{2n+1}\text{%
.}
\end{eqnarray*}
\end{theorem}

\textbf{Proof.} Consider system (\ref{equ:waszxb}) by taking $n=0,1,2,...$
as follows:%
\begin{equation*}
\begin{array}{cccc}
n=0 & \Rightarrow & x_{1}=\frac{-1}{x_{-1}y_{0}-y_{0}+1}\text{,} & y_{1}=%
\frac{-1}{y_{-1}x_{0}-x_{0}+1}\text{,} \\ 
n=1 & \Rightarrow & x_{2}=\frac{-\left( y_{-1}x_{0}-x_{0}+1\right) }{%
y_{-1}x_{0}-2x_{0}+2}\text{,} & y_{2}=\frac{-\left(
x_{-1}y_{0}-y_{0}+1\right) }{x_{-1}y_{0}-2y_{0}+2}\text{,} \\ 
n=2 & \Rightarrow & x_{3}=\frac{-\left( x_{-1}y_{0}-2y_{0}+2\right) }{%
2x_{-1}y_{0}-3y_{0}+4}\text{,} & y_{3}=\frac{-\left(
y_{-1}x_{0}-2x_{0}+2\right) }{2y_{-1}x_{0}-3x_{0}+4}\text{,} \\ 
n=3 & \Rightarrow & x_{4}=\frac{-\left( 2y_{-1}x_{0}-3x_{0}+4\right) }{%
4y_{-1}x_{0}-6x_{0}+7}\text{,} & y_{4}=\frac{-\left(
2x_{-1}y_{0}-3y_{0}+4\right) }{4x_{-1}y_{0}-6y_{0}+7}\text{,} \\ 
n=4 & \Rightarrow & x_{5}=\frac{-\left( 4x_{-1}y_{0}-y_{0}+7\right) }{%
7x_{-1}y_{0}-11y_{0}+13}\text{,} & y_{5}=\frac{-\left(
4y_{-1}x_{0}-6x_{0}+7\right) }{7y_{-1}x_{0}-11x_{0}+13}\text{,} \\ 
n=5 & \Rightarrow & x_{6}=\frac{-\left( 7y_{-1}x_{0}-11x_{0}+13\right) }{%
13y_{-1}x_{0}-20x_{0}+24}\text{,} & x_{6}=\frac{-\left(
7x_{-1}y_{0}-11y_{0}+13\right) }{13x_{-1}y_{0}-20y_{0}+24}\text{,} \\ 
& \vdots &  & 
\end{array}%
\end{equation*}%
If we keep on this process and also regard (\ref{equ:opopopop}), then the
result directly follows from a simple induction.

\begin{theorem}
\label{theorem:ewqdaxzvsftsd}The system (\ref{equ:waszxb}) has unique
negative equilibrium point $\left( \overline{x},\overline{y}\right) =\left(
d,d\right) $ and $\left( d,d\right) $ is locally asymptotically stable.
\end{theorem}

\textbf{Proof.} Clearly, equilibrium point of system (\ref{equ:waszxb}) is
the real roots of the equations%
\begin{equation}
\overline{x}=\frac{-1}{\overline{x}\left( \overline{y}-1\right) +1}\text{, \ 
}\overline{y}=\frac{-1}{\overline{y}\left( \overline{x}-1\right) +1}\text{.}
\label{equ:pppppppppppppl}
\end{equation}%
In (\ref{equ:pppppppppppppl}), after some operations, we get%
\begin{equation*}
\overline{x}=\overline{y}\text{.}
\end{equation*}%
As a result, we obtain the following equation%
\begin{equation}
\overline{x}^{3}-\overline{x}^{2}+\overline{x}+1=0\text{.}  \label{equ:ipolk}
\end{equation}%
Then, the roots of the cubic equation (\ref{equ:ipolk}) are given by%
\begin{eqnarray*}
d &=&\frac{1+\sqrt[3]{3\sqrt{33}-17}-\sqrt[3]{3\sqrt{33}+17}}{3}\text{,} \\
e &=&\frac{1+\omega \sqrt[3]{3\sqrt{33}-17}-\omega ^{2}\sqrt[3]{3\sqrt{33}+17%
}}{3}\text{,} \\
f &=&\frac{1+\omega ^{2}\sqrt[3]{3\sqrt{33}-17}-\omega \sqrt[3]{3\sqrt{33}+17%
}}{3}\text{,}
\end{eqnarray*}%
where%
\begin{equation*}
\omega =\frac{-1+i\sqrt{3}}{2}=\exp \left( 2\pi i/3\right)
\end{equation*}%
is a primitive cube root of unity.\ So, the root $d$ is only real number.
Therefore, the unique negative equilibrium point of system (\ref{equ:waszxb}%
) is $\left( \overline{x},\overline{y}\right) =\left( d,d\right) $.

Now, we show that the unique negative equilibrium point of system (\ref%
{equ:waszxb}) is locally asymptotically stable.

Let $I$ and $J$ are some intervals of real numbers.and consider the functions

\begin{equation*}
f:I^{2}\times J^{2}\rightarrow I\text{ and }g:I^{2}\times J^{2}\rightarrow J
\end{equation*}%
defined by%
\begin{equation*}
f\left( x_{n},x_{n-1},y_{n},y_{n-1}\right) =\frac{-1}{y_{n}\left(
x_{n-1}-1\right) +1},\ g\left( x_{n},x_{n-1},y_{n},y_{n-1}\right) =\frac{-1}{%
x_{n}\left( y_{n-1}-1\right) +1}\text{.}
\end{equation*}%
We consider the following transformation to build corresponding linearized
form of system (\ref{equ:waszxb})%
\begin{equation*}
\left( x_{n},x_{n-1},y_{n},y_{n-1}\right) \rightarrow \left(
f,f_{1},g,g_{1}\right) \text{,}
\end{equation*}%
where%
\begin{eqnarray*}
f\left( x_{n},x_{n-1},y_{n},y_{n-1}\right) &=&\frac{-1}{y_{n}\left(
x_{n-1}-1\right) +1}\text{,} \\
f_{1}\left( x_{n},x_{n-1},y_{n},y_{n-1}\right) &=&x_{n}\text{,} \\
g\left( x_{n},x_{n-1},y_{n},y_{n-1}\right) &=&\frac{-1}{x_{n}\left(
y_{n-1}-1\right) +1}\text{,} \\
g_{1}\left( x_{n},x_{n-1},y_{n},y_{n-1}\right) &=&y_{n}\text{.}
\end{eqnarray*}

The linearized system of system (\ref{equ:waszxb}) about the equilibrium
point $\left( d,d\right) $ under the above transformation is given as%
\begin{equation*}
X_{n+1}=BX_{n}\text{,}
\end{equation*}%
where $X_{n}=\left( x_{n},x_{n-1},y_{n},y_{n-1}\right) ^{T}$ and $B$ is a
Jacobian matrix of system (\ref{equ:waszxb}) about the equilibrium point $%
\left( d,d\right) $ and given by%
\begin{eqnarray*}
B &=&\left( 
\begin{array}{cccc}
0 & \frac{d}{\left( d\left( d-1\right) +1\right) ^{2}} & \frac{d-1}{\left(
d\left( d-1\right) +1\right) ^{2}} & 0 \\ 
1 & 0 & 0 & 0 \\ 
\frac{d-1}{\left( d\left( d-1\right) +1\right) ^{2}} & 0 & 0 & \frac{d}{%
\left( d\left( d-1\right) +1\right) ^{2}} \\ 
0 & 0 & 1 & 0%
\end{array}%
\right) \\
&=&\left( 
\begin{array}{cccc}
0 & d^{3} & -\left( 1+d\right) & 0 \\ 
1 & 0 & 0 & 0 \\ 
-\left( 1+d\right) & 0 & 0 & d^{3} \\ 
0 & 0 & 1 & 0%
\end{array}%
\right) \text{.}
\end{eqnarray*}%
Thus, we obtain the characteristic equation of the Jacobian matrix $B$ as%
\begin{equation*}
\left( d^{3}-\lambda ^{2}\right) ^{2}-\left( 1+d\right) ^{2}\lambda ^{2}=0%
\text{,}
\end{equation*}%
or%
\begin{equation*}
\left( \lambda ^{2}-\left( 1+d\right) \lambda -d^{3}\right) \left( \lambda
^{2}+\left( 1+d\right) \lambda -d^{3}\right) =0\text{.}
\end{equation*}%
Hence, it is clearly seen that numerically%
\begin{equation*}
\left \vert \lambda _{1}\right \vert =\left \vert \lambda _{2}\right \vert
=\left \vert \lambda _{3}\right \vert =\left \vert \lambda _{4}\right \vert
=0.40089<1\text{.}
\end{equation*}%
Consequently, the equilibrium point $\left( d,d\right) $ is locally
asymptotically stable.

\begin{theorem}
The equilibrium point of system (\ref{equ:waszxb}) is globally
asymptotically stable.
\end{theorem}

\textbf{Proof.} Let $\left\{ x_{n},y_{n}\right\} _{n\geq -1}$ be a solution
system (\ref{equ:waszxb}). By Theorem (\ref{theorem:ewqdaxzvsftsd}), we need
only to prove that the equilibrium point $\left( d,d\right) $ is global
attractor, that is%
\begin{equation*}
\lim_{n\rightarrow \infty }\left( x_{n},y_{n}\right) =\left( d,d\right) 
\text{.}
\end{equation*}%
From Theorem (\ref{theorem:kmuyhn}), (\ref{equ:eeeeeeeeee}) and (\ref%
{lim:wwswws}), it follows that%
\begin{eqnarray*}
\lim_{n\rightarrow \infty }x_{2n-1} &=&\lim_{n\rightarrow \infty }\frac{%
-\left( T_{2n-2}x_{-1}y_{0}+\left( T_{2n-1}-T_{2n}\right)
y_{0}+T_{2n-1}\right) }{T_{2n-1}x_{-1}y_{0}-\left( T_{2n-2}+T_{2n-1}\right)
y_{0}+T_{2n}} \\
&=&\lim_{n\rightarrow \infty }\frac{-T_{2n-2}\left( x_{-1}y_{0}+\left( \frac{%
T_{2n-1}}{T_{2n-2}}-\frac{T_{2n}}{T_{2n-2}}\right) y_{0}+\frac{T_{2n-1}}{%
T_{2n-2}}\right) }{T_{2n-1}\left( x_{-1}y_{0}-\left( \frac{T_{2n-2}}{T_{2n-1}%
}+1\right) y_{0}+\frac{T_{2n}}{T_{2n-1}}\right) } \\
&=&\left( \frac{x_{-1}y_{0}+\left( \alpha -\alpha ^{2}\right) y_{0}+\alpha }{%
x_{-1}y_{0}-\left( \frac{1}{\alpha }+1\right) y_{0}+\alpha }\right)
\lim_{n\rightarrow \infty }\frac{-T_{2n-2}}{T_{2n-1}} \\
&=&\lim_{n\rightarrow \infty }\frac{-T_{2n-2}}{T_{2n-1}} \\
&=&-\frac{1}{\alpha } \\
&=&d,
\end{eqnarray*}%
and%
\begin{eqnarray*}
\lim_{n\rightarrow \infty }x_{2n} &=&\lim_{n\rightarrow \infty }\frac{%
-\left( T_{2n-1}y_{-1}x_{0}+\left( T_{2n}-T_{2n+1}\right)
x_{0}+T_{2n}\right) }{T_{2n}y_{-1}x_{0}-\left( T_{2n-1}+T_{2n}\right)
x_{0}+T_{2n+1}} \\
&=&\lim_{n\rightarrow \infty }\frac{-T_{2n-1}\left( y_{-1}x_{0}+\left( \frac{%
T_{2n}}{T_{2n-1}}-\frac{T_{2n+1}}{T_{2n-1}}\right) x_{0}+\frac{T_{2n}}{%
T_{2n-1}}\right) }{T_{2n}\left( y_{-1}x_{0}-\left( \frac{T_{2n-1}}{T_{2n}}%
+1\right) x_{0}+\frac{T_{2n+1}}{T_{2n}}\right) } \\
&=&\left( \frac{y_{-1}x_{0}+\left( \alpha -\alpha ^{2}\right) x_{0}+\alpha }{%
y_{-1}x_{0}-\left( \frac{1}{\alpha }+1\right) x_{0}+\alpha }\right)
\lim_{n\rightarrow \infty }\frac{-T_{2n-1}}{T_{2n}} \\
&=&\lim_{n\rightarrow \infty }\frac{-T_{2n-1}}{T_{2n}} \\
&=&-\frac{1}{\alpha } \\
&=&d\text{.}
\end{eqnarray*}%
Then, we have%
\begin{equation*}
\lim_{n\rightarrow \infty }x_{n}=d\text{.}
\end{equation*}%
Similarly, we obtain%
\begin{equation*}
\lim_{n\rightarrow \infty }y_{n}=d\text{.}
\end{equation*}%
Therefore, we get%
\begin{equation*}
\lim_{n\rightarrow \infty }\left( x_{n},y_{n}\right) =\left( d,d\right) 
\text{,}
\end{equation*}%
which completes the proof.

\end{document}